\documentclass[12pt]{article}

\usepackage{amsmath,amsfonts,amssymb,amsthm,cite}
\usepackage{geometry}

\theoremstyle{plain}
\newtheorem{theorem}{Theorem}

\newtheorem{corollary}{Corollary}
\newtheorem{definition}{Definition}
\theoremstyle{definition}

\newcommand{\nonprint}[1]{}

\begin{document}

\begin{flushleft}

\textit{Conference Paper, December, 2019\\
International Workshop on the Qualitative Theory of Differential Equations\\ ''QUALITDE – 2019'', Tbilisi, Georgia}
\vspace{+0.2cm}

\textbf{Olena~Atlasiuk, Vladimir~Mikhailets} 

\small
Institute of Mathematics of the National Academy of Science of Ukraine, Kyiv, Ukraine

\textit{E-mails: \textbf{hatlasiuk@gmail.com; mikhailets@imath.kiev.ua}}

\Large

\textbf{On Linear Boundary-Value Problems for Differential Systems in Sobolev spaces}
\end{flushleft}

\normalsize

\begin{abstract}
We consider the Fredholm one-dimensional boundary-value problems in Sobolev spaces.We have obtained several important results about the indixes of functional operators, the criterion of their correct well-posedness, the criterion of the continuous dependence of the solutions of these problems on the parameter, the degree of convergence of these solutions, and sufficient constructive conditions under which the solutions of the most general class of multipoint boundary-value problems are continuous with respect to the parameter. Eeach of these boundary-value problems corresponds to a certain rectangular numerical characteristic matrix with kernel and cokernel of the same dimension as the kernel and cokernel of the boundary-value problem. The conditions for the sequence of characteristic matrices to converge are found.
\end{abstract}

Let a finite interval  $[a,b]\subset\mathbb{R}$ and parameters $\{m, \, n, \,r\} \subset \mathbb{N}, \, 1\leqslant p\leqslant \infty$, be given. By $W_p^{n}=W_p^{n}\bigl([a,b];\mathbb{C}\bigr):= \bigl\{y\in C^{n-1}[a,b] \colon y^{(n-1)}\in AC[a,b], \, y^{(n)}\in L_p[a,b]\bigr\}$ we denote a complex Sobolev space and set $W_p^0:=L_p$. This space is a Banach one with respect to the norm
$$
\bigl\|y\bigr\|_{n,p}=\sum_{k=0}^{n-1}\bigl\|y^{(k)}\bigr\|_{p}+\bigl\|y^{(n)}\bigr\|_p,
$$
where $\|\cdot\|_p$ is the norm in the space $L_p\bigl([a,b]; \mathbb{C}\bigr)$. Similarly, by $(W_p^n)^{m}:=W_p^n\bigl([a,b];\mathbb{C}^{m}\bigr)$ and $(W_p^n)^{m\times m}:=W_p^n\bigl([a,b];\mathbb{C}^{m\times m}\bigr)$
we denote Sobolev spaces of vector-valued functions and matrix-valued functions, respectively, whose ele\-ments belong to the function space $W_p^{n}$.

We consider the following linear boundary-value problem
\begin{equation}\label{1}
    Ly(t):= y'(t)+A(t)y(t)=f(t),\quad
t \in(a,b), 
\end{equation}
\begin{equation}\label{2}
    By= c, 
\end{equation}
where the matrix-valued function $A(\cdot)\in (W_{p}^{n-1})^{m\times m}$, the vector-valued function $f(\cdot)\in (W_{p}^{n-1})^{m}$, the vector $c\in\mathbb{C}^{r}$, the linear continuous operator
\begin{equation}\label{3}
B\colon(W_{p}^{n})^{m} \rightarrow\mathbb{C}^{r} 
\end{equation}
are arbitrarily chosen; and the vector-valued function $y(\cdot)\in (W_{p}^{n})^m$ is unknown.

We represent vectors and vector-valued functions in the form of columns. A solution of the boundary-value problem \eqref{1}, \eqref{2} is understood as a vector-valued function $y(\cdot)\in (W_{p}^{n})^m$ satisfying equation \eqref{1} almost everywhere on $(a,b)$ (everywhere for $n \geq 2$) and equality \eqref{2} specifying $r$ scalar boundary conditions. The solutions of equation~\eqref{1} fill the space $(W_{p}^{n})^m$ if its right-hand side $f(\cdot)$ runs through the space $(W_{p}^{n-1})^m$. Hence, the boundary condition~\eqref{2} with continuous operator~\eqref{3} is the most general condition for this equation.

It includes all known types of classical boundary conditions, namely, the Cauchy problem, two- and many-point problems, integral and mixed problems, and numerous nonclassical problems. The last class of problems may contain derivatives of the unknown functions of the~order~$k\leqslant n$.

It is known that, for $1\leq p < \infty$, every operator $B$ in \eqref{3}
admits a unique analytic representation
$$
By=\sum _{k=0}^{n-1} \alpha_{k} y^{(k)}(a)+\int_{a}^b \Phi(t)y^{(n)}(t){\rm d}t, \quad y(\cdot)\in (W_{p}^{n})^{m},
$$
where the matrices $\alpha_{k}\in\mathbb{C}^{r\times m}$ and the matrix-valued function $\Phi(\cdot)\in L_{p^{'}}\bigl([a, b]; \mathbb{C}^{r\times m}\bigr)$, \smash{$1/p + 1/p^{'}=1$}.

For $p=\infty$ this formula also defines an operator $B\in L\bigl((W_{\infty}^{n})^{m}; \mathbb{C}^{r}\bigr)$. However, there exist other operators from this class generated by the integrals over finitely additive measures.

We rewrite the inhomogeneous boundary-value problem \eqref{1}, \eqref{2} in the form of a linear operator equation $(L,B)y=(f,c)$, where $(L,B)$ is a linear operator in the pair of Banach spaces
\begin{equation}\label{4}
(L,B)\colon (W^{n}_p)^m\rightarrow (W^{n-1}_p)^m\times\mathbb{C}^r. 
\end{equation}

Recall that a linear continuous operator $T\colon X \rightarrow Y$, where $X$ and $Y$ are Banach spaces, is called a Fredholm
operator if its kernel $\ker T$ and cokernel $Y/T(X)$ are finite-dimensional. If operator $T$ is Fredholm, then its
range $T(X)$ is closed in $Y$ and the index
$$
\mathrm{ind}\,T:=\dim\ker T-\dim(Y/T(X))
$$
is finite.

\begin{theorem}\label{Th1}
The linear operator \eqref{4} is a bounded Fredholm operator with index $m-r$.
\end{theorem}

Theorem 1 allows the next refinement.

By $Y(\cdot)\in (W_p^n)^{m\times m}$ we denote a unique solution of the linear homogenous matrix equation $(LY)(t)=O_m$, $Y(a)=I_m$, where $O_m$ is the $(m \times m)$ zero matrix, and $I_m$ is the $(m \times m)$ identity matrix.

\begin{definition}\label{D1}
A rectangular numerical matrix
$
M(L,B)\in\mathbb{C}^{m \times r}
$
is characteristic for the~boundary-value problem \eqref{1}, \eqref{2} if its $j$-th column is the result of the action of the operator $B$ on the $j$-th column of $Y(\cdot)$.
\end{definition}

Here $m$ is the number of scalar differential equations of the system \eqref{1}, and $r$ is the number of scalar boundary conditions.

\begin{theorem}\label{Th2}
The dimensions of the kernel and cokernel of the operator \eqref{4} are equal to the~dimensions of the kernel and cokernel of the characteristic matrix $M(L,B)$ respectively.
\end{theorem}

Theorem \ref{Th2} implies a criterion for the invertibility of the operator \eqref{4}.

\begin{corollary}\label{C1}
The operator $(L,B)$ is invertible if and only if $r=m$ and the matrix $M(L,B)$ is nondegenerate.
\end{corollary}

Let us consider parameterized by number $\varepsilon \in [0,\varepsilon_0)$, $\varepsilon_0>0$, linear boundary-value problem
\begin{equation}\label{5}
 L(\varepsilon)y(t;\varepsilon):= y'(t;\varepsilon) + A(t;\varepsilon) y(t;\varepsilon) = f(t;\varepsilon), \quad t\in (a,b),
\end{equation}
\begin{equation}\label{6} 
 B(\varepsilon)y(\cdot;\varepsilon) = c(\varepsilon),
\end{equation}
where for every fixed $\varepsilon$ the matrix-valued function $A(\cdot;\varepsilon) \in (W^{n-1}_p)^{m\times m}$, the vector-valued function \smash{$f(\cdot;\varepsilon) \in (W^{n-1}_p)^m$}, the vector $c(\varepsilon)\in\mathbb{C}^m$, $B(\varepsilon)$~is the linear continuous operator $B(\varepsilon) \colon (W^{n}_p)^m\rightarrow \mathbb{C}^m$, and the solution (the unknown vector-valued function) $y(\cdot;\varepsilon) \in (W^{n}_p)^m$.

It follows from Theorem \ref{Th2} that the boundary-value problem \eqref{5}, \eqref{6} is Fredholm with index zero.

\begin{definition}\label{D2}
A solution of the boundary-value problem \eqref{5}, \eqref{6} continuously depends on the parameter $\varepsilon$ for $\varepsilon=0$ if the following conditions are satisfied:
\begin{itemize}
\item [$(\ast)$] there exists a positive number $\varepsilon_{1}<\varepsilon_{0}$ such that, for any $\varepsilon\in[0,\varepsilon_{1})$ and arbitrary right-hand sides $f(\cdot;\varepsilon)\in (W^{n-1}_p)^{m}$ and  $c(\varepsilon)\in\mathbb{C}^{m}$ this problem has a unique solution $y(\cdot;\varepsilon)$ that belongs to the space $(W^{n}_p)^{m}$;

\item [$(\ast\ast)$] the convergence of the right-hand sides $f(\cdot;\varepsilon)\to f(\cdot;0)$ in $(W_p^{n-1})^{m}$ and $c(\varepsilon)\to c(0)$ in $\mathbb{C}^{m}$ as $\varepsilon\to0+$ implies the convergence of the solutions $y(\cdot;\varepsilon)\to y(\cdot;0)$ in $(W^{n}_p)^{m}$.
\end{itemize}
\end{definition}

Consider the following conditions as $\varepsilon\to0+$:
\begin{itemize}
  \item [(0)] limiting homogeneous boundary-value problem
$$L(0)y(t,0)=0,\quad t\in(a, b),\quad B(0)y(\cdot,0)=0$$
has only the trivial solution;
  \item [(I)] $A(\cdot,\varepsilon)\to A(\cdot,0)$ in the space $(W^{n-1}_{p})^{m\times m}$;
  \item [(II)] $B(\varepsilon)y\to B(0)y$ in $\mathbb{C}^{m}$ for any $y\in(W^{n}_{p})^m$.
\end{itemize}
\begin{theorem}\label{Th3}
A solution of the boundary-value problem \eqref{5}, \eqref{6} continuously depends on the~parameter $\varepsilon$ for $\varepsilon=0$ if and only if it satisfies condition (0) and the conditions (I) and (II).
\end{theorem}
Consider the following quantities:
\begin{equation}\label{7}
\bigl\|y(\cdot;0)-y(\cdot;\varepsilon)\bigr\|_{n,p}, 
\end{equation}
\begin{equation}\label{8}
\widetilde{d}_{n-1,p}(\varepsilon):=
\bigl\|L(\varepsilon)y(\cdot;0)-f(\cdot;\varepsilon)\bigr\|_{n-1,p}+
\bigl\|B(\varepsilon)y(\cdot;0)-c(\varepsilon)\bigr\|_{\mathbb{C}^{m}}, 
\end{equation}
where \eqref{7} is the error and \eqref{8} is the discrepancy of the solution $y(\cdot;\varepsilon)$ of the boundary-value problem~\eqref{5}, \eqref{6} if $y(\cdot;\varepsilon)$ is its exact solution and $y(\cdot;0)$ is an approximate solution of the~problem.

\begin{theorem}\label{Th4}
Suppose that the boundary-value problem \eqref{5}, \eqref{6} satisfies conditions (0), (I) and (II). Then there exist the positive quantities $\varepsilon_{2}<\varepsilon_{1}$ and $\gamma_{1}$, $\gamma_{2}$ such that, for any  $\varepsilon\in(0,\varepsilon_{2})$, the following two-sided estimate is true:
$$
\gamma_{1}\,\widetilde{d}_{n-1,p}(\varepsilon)
\leq\bigl\|y(\cdot;0)-y(\cdot;\varepsilon)\bigr\|_{n,p}\leq
\gamma_{2}\,\widetilde{d}_{n-1,p}(\varepsilon),
$$
where the quantities $\varepsilon_{2}$, $\gamma_{1}$, and $\gamma_{2}$ do not depend of $y(\cdot;\varepsilon)$ and $y(\cdot;0)$.
\end{theorem}

Acording to this theorem, the error and discrepancy of the solution $y(\cdot;\varepsilon)$ of the boundary-value problem \eqref{5}, \eqref{6} have the same order of smallness.

For any $\varepsilon \in [0,\varepsilon_0)$, $\varepsilon_0>0$, we associate with the system~\eqref{5} multi-point Fredholm boundary condition
\begin{equation}\label{9}
B(\varepsilon)y(\cdot,\varepsilon)= \sum\limits_{j=0}^{r}\sum\limits_{k=1}^{\omega_j(\varepsilon)}
\sum\limits_{l=0}^{n}{\beta_{j,k}^{(l)}
(\varepsilon)y^{(l)}(t_{j,k}(\varepsilon),\varepsilon)}=q(\varepsilon), 
\end{equation}
where the numbers $\{r, \omega_j(\varepsilon)\} \subset \mathbb{N}$, vectors $q(\varepsilon)\in\mathbb{C}^{m}$,  matrices $\beta_{j,k}^{(l)}(\varepsilon)\in\mathbb{C}^{m\times m}$, and points $\{t_{j}, t_{j,k}(\varepsilon)\} \subset[a,b]$ are arbitrarily given.

It is not assumed that the coefficients $A(\cdot,\varepsilon)$, $\beta_{j,k}^{(l)}(\varepsilon)$ or points $t_{j,k}(\varepsilon)$ have a certain regularity on the parameter~$\varepsilon$ as $\varepsilon >0$. It will be required that for each fixed $j\in\{1,\ldots,r\}$ all the points $t_{j,k}(\varepsilon)$ have a common limit as $\varepsilon\to0+$, but for the zero-point series $t_{0,k}(\varepsilon)$ this requirement will not be necessary.

The solution $y=y(\cdot,\varepsilon)$ of the multi-point boundary-value problem \eqref{5}, \eqref{9} is continuous on the parameter $\varepsilon$ if it exists, is unique, and satisfies the limit relation
\begin{equation}\label{10}
\bigl\|y(\cdot,\varepsilon)-y(\cdot,0) \bigr\|_{n,p} \to 0 \quad \mbox{as} \quad \varepsilon\to0+. 
\end{equation}

Consider the following assumptions as $\varepsilon\to0+$ and $p=\infty$:
\begin{itemize}
    \item [$(\alpha)$] $t_{j,k}(\varepsilon)\to t_{j}$ for all $j\in\{1,\ldots,r\}$, and $k\in\{1,\ldots,\omega_j(\varepsilon)\};$

    \item [$(\beta)$] $\sum\limits_{k=1}^{\omega_j(\varepsilon)}
\beta_{j,k}^{(l)}(\varepsilon)\to\beta_{j}^{(l)}$ for all $j\in\{1,\ldots,r\}$, and $l\in\{0,\ldots,n\}$;

    \item [$(\gamma)$] $\sum\limits_{k=1}^{\omega_j(\varepsilon)}\bigl\|\beta_{j,k}^{(l)}(\varepsilon)\bigr\|
\bigl|t_{j,k}(\varepsilon)-t_j\bigr|\to0$ for all $j\in\{1,\ldots,r\}$, $k\in\{1,\ldots,\omega_j(\varepsilon)\}$, and $l\in\{0,\ldots,n\}$;

      \item [$(\delta)$] $\sum\limits_{k=1}^{\omega_0(\varepsilon)}\bigl\|\beta_{0,k}^{(l)}(\varepsilon)\bigr\|\to0$ for all $k\in\{1,\ldots,\omega_0(\varepsilon)\}$, and $l\in\{0,\ldots,n\}$.
\end{itemize}

Assumptions $(\beta)$ and $(\gamma)$ imply that the norms of the coefficients $\beta_{j,k}^{(l)}(\varepsilon)$ can increase as $\varepsilon\to0+$, but not too fast.

\begin{theorem}\label{Th5}
Let the boundary-value problem \eqref{5}, \eqref{9} for $p= \infty$ satisfy the assumptions $(\alpha)$, $(\beta)$, $(\gamma)$, $(\delta)$. Then it satisfies the limit condition~(II). If, moreover, the conditions (0) and (I) are fulfilled, then for a sufficiently small $\varepsilon$ its solution exists, is unique and satisfies the limit relation \eqref{10}.
\end{theorem}

Consider also the following assumptions as $\varepsilon\to0+$ and $1\leqslant p< \infty$:
\begin{itemize}
       \item[$(\gamma_p)$] $\sum\limits_{k=1}^{\omega_j(\varepsilon)}\bigl\|\beta_{j,k}^{(n)}(\varepsilon)\bigr\|
\bigl|t_{j,k}(\varepsilon)-t_j\bigr|^{1/p^{'}}=O(1)$ for all $j\in\{1,\ldots,r\}$, and $k\in\{1,\ldots,\omega_j(\varepsilon)\}$;

    \item[$(\gamma')$] $\sum\limits_{k=1}^{\omega_j(\varepsilon)}\bigl\|\beta_{j,k}^{(l)}(\varepsilon)\bigr\|
\bigl|t_{j,k}(\varepsilon)-t_j\bigr|\to0$ for all $j\in\{1,\ldots,r\}$, $k\in\{1,\ldots,\omega_j(\varepsilon)\}$, and
 \mbox{$l\in\{0,\ldots,n-1\}$}.
    \end{itemize}

\begin{theorem}\label{Th6} 
Let the boundary-value problem \eqref{5}, \eqref{9} for $1\leqslant p< \infty$ satisfy the assumptions $(\alpha)$, $(\beta)$, $(\gamma_p)$, $(\gamma')$, $(\delta)$. Then it satisfies the limit condition~(II). If, moreover, the conditions (0) and (I) are fulfilled, then for a sufficiently small $\varepsilon$ its solution exists, is unique and satisfies the limit relation \eqref{10}.
\end{theorem}

The results are published in~\cite{AtlMikh20181, AtlMikh20182, AtlMikh20183, AtlMikh20184}. They allow extension for the systems of  differential equations of higher order \cite{GKM2017} and for boundary-value problems in H\"older spaces \cite{MMS2016}.

\end{document}